\documentclass [10pt,reqno]{amsart}
\usepackage {hannah}
\usepackage {color,graphicx}
\usepackage {multirow}

\newcommand {\brk}{\rule {1mm}{0mm}}
\newcommand {\df}[1]{\emph {#1}}

\newcommand {\calM}{{\mathcal M}}
\newcommand {\calP}{{\mathcal P}}
\newcommand {\PP}{{\mathbb P}}
\newcommand {\RR}{{\mathbb R}}
\newcommand {\ZZ}{{\mathbb Z}}

\DeclareMathOperator {\ev}{ev}
\DeclareMathOperator {\ft}{f\/t}
\DeclareMathOperator {\val}{val}

\newcommand {\dunion}{\,\mbox {\raisebox{0.25ex}{$\cdot$} \kern-1.83ex $\cup$}
  \,}

\title [The WDVV equations in tropical geometry]{Kontsevich's formula and the
  WDVV equations in tropical geometry}
\author {Andreas Gathmann and Hannah Markwig}
\address {Andreas Gathmann, Fachbereich Mathematik, Technische Universit\"at
  Kaiserslautern, Postfach 3049, 67653 Kaiserslautern, Germany}
\email {andreas@mathematik.uni-kl.de}
\address {Hannah Markwig, Fachbereich Mathematik, Technische Universit\"at
  Kaiserslautern, Postfach 3049, 67653 Kaiserslautern, Germany}
\email {markwig@mathematik.uni-kl.de}
\thanks {\emph {2000 Mathematics Subject Classification:} Primary 14N35, 51M20,
  Secondary 14N10}
\thanks {The second author has been funded by the DFG grant Ga 636/2.}
\keywords {Tropical geometry, enumerative geometry, Gromov-Witten theory}

\begin {document}

\begin {abstract}
  Using Gromov-Witten theory the numbers of complex plane rational curves of
  degree $d$ through $ 3d-1 $ general given points can be computed recursively
  with Kontsevich's formula that follows from the so-called WDVV equations. In
  this paper we establish the same results entirely in the language of tropical
  geometry. In particular this shows how the concepts of moduli spaces of
  stable curves and maps, (evaluation and forgetful) morphisms, intersection
  multiplicities and their invariance under deformations can be carried over to
  the tropical world.
\end {abstract}

\maketitle


\section {Introduction}

For $ d \ge 1 $ let $ N_d $ be the number of rational curves in the complex
projective plane $ \PP^2 $ that pass through $ 3d-1 $ given points in general
position. About 10 years ago Kontsevich has shown that these numbers are given
recursively by the initial value $ N_1 = 1 $ and the equation
  \[ N_d = \sum_{\substack {d_1+d_2=d \\ d_1,d_2>0}} \left(
       d_1^2 d_2^2 \binom {3d-4}{3d_1-2} - d_1^3 d_2 \binom {3d-4}{3d_1-1}
     \right) N_{d_1} N_{d_2} \]
for $ d>1 $ (see \cite {KM94} claim 5.2.1). The main tool in deriving this
formula is the so-called WDVV equations, i.e.\ the associativity equations of
quantum cohomology. Stated in modern terms the idea of these equations is as
follows: plane rational curves of degree $d$ are parametrized by the moduli
spaces of stable maps $ \bar M_{0,n} (\PP^2,d) $ whose points are in bijection
to tuples $ (C,x_1,\dots,x_n,f) $ where $ x_1,\dots,x_n $ are distinct smooth
points on a rational nodal curve $C$ and $ f:C \to \PP^2 $ is a morphism of
degree $d$ (with a stability condition). If $ n \ge 4 $ there is a ``forgetful
map'' $ \pi: \bar M_{0,n} (\PP^2,d) \to \bar M_{0,4} $ that sends a stable map
$ (C,x_1,\dots,x_n,f) $ to (the stabilization of) $ (C,x_1,\dots,x_4) $. The
important point is now that the moduli space $ \bar M_{0,4} $ of 4-pointed
rational stable curves is simply a projective line. Therefore the two points

\begin {center} \input {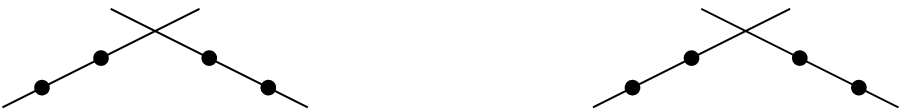} \end {center}

of $ \bar M_{0,4} $ are linearly equivalent divisors, and hence so are their
inverse images $ D_{12|34} $ and $ D_{13|24} $ under $ \pi $. The divisor $
D_{12|34} $ in $ \bar M_{0,n}(\PP^2,d) $ (and similarly of course $ D_{13|24}
$) can be described explicitly as the locus of all reducible stable maps with
two components such that the marked points $ x_1,x_2 $ lie on one component and
$ x_3,x_4 $ on the other. It is of course reducible since there are many
combinatorial choices for such curves: the degree and the remaining marked
points can be distributed onto the two components in an arbitrary way.

All that remains to be done now is to intersect the equation $ [D_{12|34}] =
[D_{13|24}] $ of divisor classes with cycles of dimension 1 in $ \bar M_{0,n}
(\PP^2,d) $ to get some equations between numbers. Specifically, to get
Kontsevich's formula one chooses $ n=3d $ and intersects the above divisors
with the conditions that the stable maps pass through two given lines at $ x_1
$ and $ x_2 $ and through given points in $ \PP^2 $ at all other $ x_i $. The
resulting equation can be seen to be precisely the recursion formula stated at
the beginning of the introduction: the sum corresponds to the possible
splittings of the degree of the curves onto their two components, the binomial
coefficients correspond to the distribution of the marked points $ x_i $ with $
i>4 $, and the various factors of $ d_1 $ and $ d_2 $ correspond to the
intersection points of the two components with each other and with the two
chosen lines (for more details see e.g.\ \cite {CK99} section 7.4.2).

The goal of this paper is to establish the same results in tropical geometry.
In contrast to most enumerative applications of tropical geometry known so far
it is absolutely crucial for this to work that we pick the ``correct''
definition of (moduli spaces of) tropical curves even for somewhat degenerated
curves.

To describe our definition let us start with abstract tropical curves, i.e.\
curves that are not embedded in some ambient space. An abstract tropical curve
is simply an abstract connected graph $ \Gamma $ obtained by glueing closed
(not necessarily bounded) real intervals together at their boundary points in
such a way that every vertex has valence at least 3. In particular, every
bounded edge of such an abstract tropical curve has an intrinsic length.
Following an idea of Mikhalkin \cite {Mi06} the unbounded ends of $ \Gamma $
will be labeled and called the marked points of the curve. The most important
example for our applications is the following:

\begin {example} \label {ex-m4}
  A 4-marked rational tropical curve (i.e.\ an element of the tropical
  analogue of $ \bar M_{0,4} $ that we will denote by $ \calM_4 $) is simply a
  tree graph with 4 unbounded ends. There are four possible combinatorial types
  for this:

  \begin {center} \input {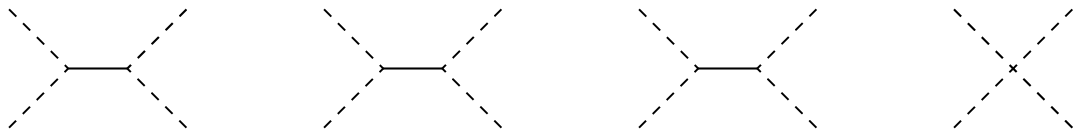} \end {center}

  (In this paper we will always draw the unbounded ends corresponding to marked
  points as dotted lines.) In the types (A) to (C) the bounded edge has an
  intrinsic length $l$; so each of these types leads to a stratum of $ \calM_4
  $ isomorphic to $ \RR_{>0} $ parametrized by this length. The last type (D)
  is simply a point in $ \calM_4 $ that can be seen as the boundary point in $
  \calM_4 $ where the other three strata meet. Therefore $ \calM_4 $ can be
  thought of as three unbounded rays meeting in a point --- note that this is
  again a rational tropical curve!

  \begin {center} \input {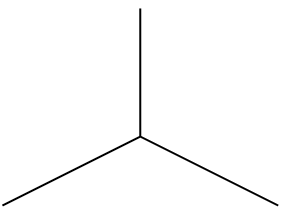} \end {center}
\end {example}

Let us now move on to plane tropical curves. As in the complex case we will
adopt the ``stable map picture'' and consider maps from an abstract tropical
curve to $ \RR^2 $ rather than embedded tropical curves. More precisely, an
$n$-marked plane tropical curve will be a tuple $ (\Gamma,x_1,\dots,x_n,h) $,
where $ \Gamma $ is an abstract tropical curve, $ x_1,\dots,x_n $ are distinct
unbounded ends of $ \Gamma $, and $ h: \Gamma \to \RR^2 $ is a continuous map
such that
\begin {enumerate}
\item on each edge of $ \Gamma $ the map $h$ is of the form $ h(t) = a + t
  \cdot v $ for some $ a \in \RR^2 $ and $ v \in \ZZ^2 $ (``$h$ is affine
  linear with integer direction vector $v$'');
\item for each vertex $V$ of $ \Gamma $ the direction vectors of the
  edges around $V$ sum up to zero (the ``balancing condition'');
\item the direction vectors of all unbounded edges corresponding to the marked
  points are zero (``every marked point is contracted to a point in $ \RR^2 $
  by $h$'').
\end {enumerate}
Note that it is explicitly allowed that $h$ contracts an edge $E$ of $ \Gamma $
to a point. If this is the case and $E$ is a bounded edge then the intrinsic
length of $E$ can vary arbitrarily without changing the image curve $ h(\Gamma)
$. This is of course the feature of ``moduli in contracted components'' that we
know well from the ordinary complex moduli spaces of stable maps.

\begin {example} \label {ex-stable-map}
  The following picture shows an example of a 4-marked plane tropical curve of
  degree 2, i.e.\ of an element of the tropical analogue of $ \bar M_{0,4}
  (\PP^2,2) $ that we will denote by $ \calM_{2,4} $. Note that at each marked
  point the balancing condition ensures that the two other edges meeting at the
  corresponding vertex are mapped to the same line in $ \RR^2 $.

  \begin {center} \input {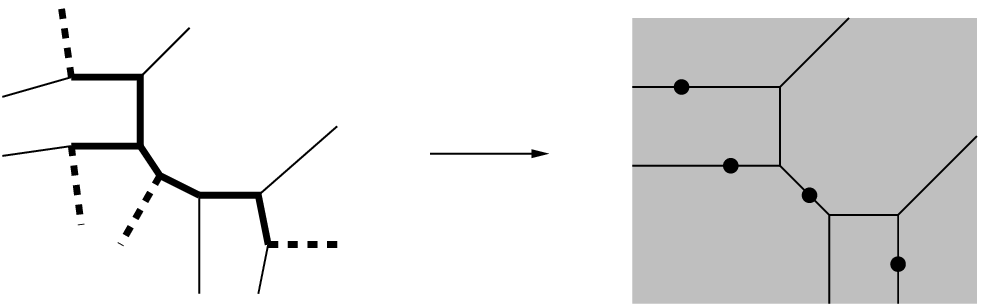} \end {center}

  It is easy to see from this picture already that the tropical moduli spaces $
  \calM_{d,n} $ of plane curves of degree $d$ with $ n \ge 4 $ marked points
  admit forgetful maps to $ \calM_4 $: given an $n$-marked plane tropical
  curve $ (\Gamma,x_1,\dots,x_n,h) $ we simply forget the map $h$, take the
  minimal connected subgraph of $ \Gamma $ that contains $ x_1,\dots,x_4 $, and 
  ``straighten'' this graph to obtain an element of $ \calM_4 $. In the picture
  above we simply obtain the ``straightened version'' of the subgraph drawn in
  bold, i.e.\ the element of $ \calM_4 $ of type (A) (in the notation of
  example \ref {ex-m4}) with length parameter $l$ as indicated in the picture.
\end {example}

The next thing we would like to do is to say that the inverse images of two
points in $ \calM_4 $ under this forgetful map are ``linearly equivalent
divisors''. However, there is unfortunately no theory of divisors in tropical
geometry yet. To solve this problem we will first impose all incidence
conditions as needed for Kontsevich's formula and then only prove that the
(suitably weighted) number of plane tropical curves satisfying all these
conditions and mapping to a given point in $ \calM_4 $ does not depend on
this choice of point. The idea to prove this is precisely the same as for the
independence of the incidence conditions in \cite {GM051} (although the
multiplicity with which the curves have to be counted has to be adapted to the
new situation).

We will then apply this result to the two curves in $ \calM_4 $ that are
of type (A) resp.\ (B) above and have a fixed very large length parameter $l$.
We will see that such very large lengths in $ \calM_4 $ can only occur if
there is a contracted bounded edge (of a very large length) somewhere as in the
following example:

\begin {example} \label {ex-conic-red}
  Let $C$ be a plane tropical curve with a bounded contracted edge $E$.

  \begin {center} \input {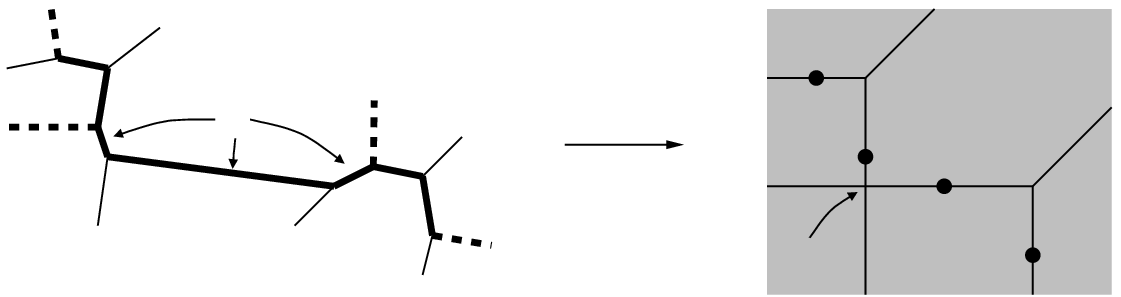} \end {center}

  In this picture the parameter $l$ is the sum of the intrinsic lengths of the
  three marked edges, in particular it is very large if the intrinsic length of
  $E$ is. By the balancing condition it follows that locally around $ P=h(E) $
  the tropical curve must be a union of two lines through $P$, i.e.\ that the
  tropical curve becomes ``reducible'' with two components meeting in $P$ (in
  the picture above we have a union of two tropical lines).
\end {example}

Hence we get the same types of splitting of the curves into two components as
in the complex picture --- and thus the same resulting formula for the
(tropical) numbers $ N_d $.

Our result shows once again quite clearly that it is possible to carry many
concepts from classical complex geometry over to the tropical world: moduli
spaces of curves and stable maps, morphisms, divisors and divisor classes,
intersection multiplicities, and so on. Even if we only make these
constructions in the specific cases needed for Kontsevich's formula we hope
that our paper will be useful to find the correct definitions of these concepts
in the general tropical setting. It should also be quite easy to generalize
our results to other cases, e.g.\ to tropical curves of other degrees
(corresponding to complex curves in toric surfaces) or in higher-dimensional
spaces. Work in this direction is in progress.

This paper is organized as follows: in section \ref {tropcurves} we define the
moduli spaces of abstract and plane tropical curves that we will work with
later. They have the structure of (finite) polyhedral complexes. For morphisms
between such complexes we then define the concepts of multiplicity and degree
in section \ref {tropmult}. We show that these notions specialize to
Mikhalkin's well-known ``multiplicities of plane tropical curves'' when applied
to the evaluation maps. In section \ref {forget} we apply the same techniques
to the forgetful maps described above. In particular, we show that the numbers
of tropical curves satisfying given incidence conditions and mapping to a given
point in $ \calM_4 $ do not depend on this choice of point in $ \calM_4 $.
Finally, we apply this result to two different points in $ \calM_4 $ to derive
Kontsevich's formula in section \ref {kontsevich}.


\section {Abstract and plane tropical curves} \label {tropcurves}

In this section we will mainly define the moduli spaces of (abstract and plane)
tropical curves that we will work with later. Our definitions here differ
slightly from our earlier ones in \cite {GM051}. A common feature of both
definitions is that we will always consider a plane curve to be a
``parametrized tropical curve'', i.e.\ a graph $ \Gamma $ with a map $h$ to the
plane rather than an embedded tropical curve. In contrast to our earlier work
however it is now explicitly allowed (and crucial for our arguments to work)
that the map $h$ contracts some edges of $ \Gamma $ to a point. Moreover,
following Mikhalkin \cite {Mi06} marked points will be contracted unbounded
ends instead of just markings. For simplicity we will only give the definitions
here for rational curves.

\begin {definition}[Graphs] \label {def-graph} \brk

  \vspace*{-3ex}

  \begin {enumerate}
  \item \label {def-graph-a}
    Let $ I_1,\dots,I_n \subset \RR $ be a finite set of closed, bounded or
    half-bounded real intervals. We pick some (not necessarily distinct)
    boundary points $ P_1,\dots,P_k,Q_1,\dots,Q_k \in I_1 \dunion \cdots
    \dunion I_n $ of these intervals. The topological space $ \Gamma $ obtained
    by identifying $ P_i $ with $ Q_i $ for all $ i=1,\dots,k $ in $ I_1
    \dunion \cdots \dunion I_n $ is called a \df {graph}. As usual, the \df
    {genus} of $ \Gamma $ is simply its first Betti number $ \dim H_1
    (\Gamma,\RR) $.
  \item \label {def-graph-b}
    For a graph $ \Gamma $ the boundary points of the intervals $ I_1,\dots,
    I_n $ are called the \df {flags}, their image points in $ \Gamma $ the \df
    {vertices} of $ \Gamma $. If $F$ is such a flag then its image vertex in $
    \Gamma $ will be denoted $ \partial F $. For a vertex $V$ the number of
    flags $F$ with $ \partial F = V $ is called the \df {valence} of $V$ and
    denoted $ \val V $. We denote by $ \Gamma^0 $ and $ \Gamma' $ the sets of
    vertices and flags of $ \Gamma $, respectively.
  \item \label {def-graph-c}
    The open intervals $ I_1^\circ,\dots,I_n^\circ $ are naturally open subsets
    of $ \Gamma $; they are called the \df {edges} of $ \Gamma $. An edge will
    be called bounded (resp.\ unbounded) if its corresponding open interval is.
    We denote by $ \Gamma^1 $ (resp.\ $ \Gamma^1_0 $ and $ \Gamma^1_\infty $)
    the set of edges (resp.\ bounded and unbounded edges) of $ \Gamma $. Every
    flag $ F \in \Gamma' $ belongs to exactly one edge that we will denote by $
    [F] \in \Gamma^1 $. The unbounded edges will also be called the \df {ends}
    of $ \Gamma $.
  \end {enumerate}
\end {definition}

\begin {definition}[Abstract tropical curves] \label {def-tropcurve}
  A (rational, abstract) tropical curve is a connected graph $ \Gamma $ of
  genus 0 all of whose vertices have valence at least 3. An \df {$n$-marked
  tropical curve} is a tuple $ (\Gamma,x_1,\dots,x_n) $ where $ \Gamma $ is a
  tropical curve and $ x_1,\dots,x_n \in \Gamma^1_\infty $ are distinct
  unbounded edges of $ \Gamma $. Two such marked tropical curves $ (\Gamma,
  x_1,\dots,x_n) $ and $ (\tilde \Gamma, \tilde x_1,\dots,\tilde x_n) $ are
  called isomorphic (and will from now on be identified) if there is a
  homeomorphism $ \Gamma \to \tilde \Gamma $ mapping $ x_i $ to $ \tilde x_i $
  for all $i$ and such that every edge of $ \Gamma $ is mapped bijectively onto
  an edge of $ \tilde \Gamma $ by an affine map of slope $ \pm 1 $, i.e.\ by a
  map of the form $ t \mapsto a \pm t $ for some $ a \in \RR $. The space of
  all $n$-marked tropical curves (modulo isomorphisms) with precisely $n$
  unbounded edges will be denoted $ \calM_n $. (It can be thought of as a
  tropical analogue of the moduli space $ \bar M_{0,n} $ of $n$-pointed stable
  rational curves.)
\end {definition}

\begin {example} \label {m03}
  We have $ \calM_n = \emptyset $ for $ n<3 $ since any graph of genus 0 all of
  whose vertices have valence at least 3 must have at least 3 unbounded edges.
  For $ n=3 $ unbounded edges there is exactly one such tropical curve, namely

  \begin {center} \input {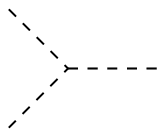} \end {center}

  (in this paper we will always draw the unbounded edges corresponding to the
  markings $ x_i $ as dotted lines). Hence $ \calM_3 $ is simply a point.
\end {example}

\begin {remark} \label {can-par}
  The isomorphism condition of definition \ref {def-tropcurve} means that every
  edge of a marked tropical curve has a parametrization as an interval in $ \RR
  $ that is unique up to translations and sign. In particular, every bounded
  edge $E$ of a tropical curve has an intrinsic \df {length} that we will
  denote by $ l(E) \in \RR_{>0} $.

  One way to fix this translation and sign ambiguity is to pick a flag $F$ of
  the edge $E$: there is then a unique choice of parametrization such that the
  corresponding closed interval is $ [0,l(E)] $ (or $ [0,\infty) $ for
  unbounded edges), with the chosen flag $F$ being the zero point of this
  interval. We will call this the \df {canonical parametrization} of $E$ with
  respect to the flag $F$.
\end {remark}

\begin {example} \label {m04}
  The moduli space $ \calM_4 $ is simply a rational tropical curve with 3
  ends --- see example \ref {ex-m4}.
\end {example}

\begin {definition}[Plane tropical curves] \label {def-planecurve} \brk

  \vspace*{-3ex}

  \begin {enumerate} \parindent 0mm \parskip 0.3ex
  \item \label {def-planecurve-a}
    Let $ n \ge 0 $ be an integer. An \df {$n$-marked plane tropical curve} is
    a tuple $ (\Gamma,x_1,\dots,x_n,h) $, where $ \Gamma $ is an abstract
    tropical curve, $ x_1,\dots,x_n \in \Gamma^1_\infty $ are distinct
    unbounded edges of $ \Gamma $, and $ h: \Gamma \to \RR^2 $ is a continuous
    map, such that:
    \begin {enumerate}
    \item [(i)] On each edge of $ \Gamma $ the map $h$ is of the form $ h(t) =
      a + t \cdot v $ for some $ a \in \RR^2 $ and $ v \in \ZZ^2 $ (i.e.\ ``$h$
      is affine linear with rational slope''). The integral vector $v$
      occurring in this equation if we pick for $E$ the canonical
      parametrization with respect to a chosen flag $F$ of $E$ (see remark \ref
      {can-par}) will be denoted $ v(F) $ and called the \df {direction} of
      $F$.
    \item [(ii)] For every vertex $V$ of $ \Gamma $ we have the \df {balancing
      condition}
        \[ \sum_{F \in \Gamma': \partial F = V} v(F) = 0. \]
    \item [(iii)] Each of the unbounded edges $ x_1,\dots,x_n \in
      \Gamma^1_\infty $ is mapped to a point in $ \RR^2 $ by $h$ (i.e.\ $
      v(F)=0 $ for the corresponding flags).
    \end {enumerate}
  \item \label {def-planecurve-b}
    Two $n$-marked plane tropical curves $ (\Gamma,x_1,\dots,x_n,h) $ and
    $ (\tilde \Gamma, \tilde x_1,\dots,\tilde x_n,\tilde h) $ are called
    isomorphic (and will from now on be identified) if there is an isomorphism
    $ \varphi: (\Gamma,x_1,\dots,x_n) \to (\tilde \Gamma,\tilde x_1,\dots,
    \tilde x_n) $ of the underlying abstract curves as in definition \ref
    {def-tropcurve} such that $ \tilde h \circ \varphi = h $.
  \item \label {def-planecurve-c}
    The \df {degree} of an $n$-marked plane tropical curve is defined to be the
    multiset $ \Delta = \{ v(F);\; [F] \in \Gamma^1_\infty \backslash \{
    x_1,\dots,x_n \} \} $ of directions of its non-marked unbounded edges. If
    this degree consists of the vectors $ (-1,0) $, $ (0,-1) $, $ (1,1) $ each
    $d$ times then we simply say that the degree of the curve is $d$. The space
    of all $n$-marked plane tropical curves of degree $ \Delta $ (resp.\ $d$)
    will be denoted $ \calM_{\Delta,n} $ (resp.\ $ \calM_{d,n} $). It can be
    thought of as a tropical analogue of the moduli spaces of stable maps to
    toric surfaces (resp.\ the projective plane).
  \end {enumerate}
\end {definition}

\begin {remark}
  For a concrete example of a marked plane tropical curve see example \ref
  {ex-stable-map}.

  Note that the map $h$ of a marked plane tropical curve $
  (\Gamma,x_1,\dots,x_n,h) $ need not be injective on the edges of $ \Gamma $:
  it may happen that $ v(F)= 0 $ for a flag $F$, i.e.\ that the corresponding
  edge is contracted to a point. Of course it follows then in such a case that
  the remaining flags around the vertex $ \partial F $ satisfy the balancing
  condition themselves. If $ \partial F $ is a 3-valent vertex this means that
  the other two flags around this vertex are negatives of each other, i.e.\
  that the image $ h(\Gamma) $ in $ \RR^2 $ is just a straight line locally
  around this vertex.

  This applies in particular to the marked unbounded edges $ x_1,\dots,x_n $ as
  they are required to be contracted by $h$. They can therefore be seen as
  tropical analogues of marked points in the ordinary complex moduli spaces of
  stable maps. By abuse of notation we will therefore often refer to these
  marked unbounded edges as ``marked points'' in the rest of the paper.

  Note also that contracted bounded edges lead to ``hidden moduli parameters''
  of plane tropical curves: if we vary the length of a contracted bounded edge
  then we arrive at a continuous family of different plane tropical curves
  whose images in $ \RR^2 $ are all the same. This feature of moduli in
  contracted components is of course well-known from the complex moduli spaces
  of stable maps.
\end {remark}

\begin {remark}
  If the direction $ v(F) \in \ZZ^2 $ of a flag $F$ of a plane tropical curve
  is not equal to zero then it can be written uniquely as a positive integer
  times a primitive integral vector. This positive integer is what is usually
  called the \df {weight} of the corresponding edge. In this paper we will not
  use this notation however since it seems more natural for our applications
  not to split up the direction vectors in this way.
\end {remark}

The following results about the structure of the spaces $ \calM_n $ and $
\calM_{\Delta,n} $ are very similar to those in \cite {GM051}, albeit much
simpler.

\begin {definition}[Combinatorial types]
  The \df {combinatorial type} of a marked tropical curve $ (\Gamma,x_1,\dots,
  x_n) $ is defined to be the homeomorphism class of $ \Gamma $ relative $
  x_1,\dots,x_n $ (i.e.\ the data of $ (\Gamma,x_1,\dots,x_n) $ modulo
  homeomorphisms of $ \Gamma $ that map each $ x_i $ to itself). The
  combinatorial type of a marked plane tropical curve $
  (\Gamma,x_1,\dots,x_n,h) $ is the data of the combinatorial type of the
  marked tropical curve $ (\Gamma,x_1,\dots,x_n) $ together with the direction
  vectors $ v(F) $ for all flags $ F \in \Gamma' $. In both cases the
  codimension of such a type $ \alpha $ is defined to be
    \[ \codim \alpha := \sum_{V\in \Gamma^0} (\val V-3). \]
  We denote by $ \calM_n^\alpha $ (resp.\ $ \calM_{\Delta,n}^{\alpha} $) the
  subset of $ \calM_n $ (resp.\ $ \calM_{\Delta,n} $) that corresponds to
  marked tropical curves of type $\alpha$.
\end {definition}

\begin {lemma} \label {fintypes}
  For all $n$ and $ \Delta $ there are only finitely many combinatorial types
  occurring in the spaces $ \calM_n $ and $ \calM_{\Delta,n} $.
\end {lemma}

\begin {proof}
  The statement is obvious for $ \calM_n $. For $ \calM_{\Delta,n} $ we just
  note in addition that by \cite {Mi03} proposition 3.11 the image $ h(\Gamma)
  $ is dual to a lattice subdivision of the polygon associated to $ \Delta $.
  In particular, this means that the absolute value of the entries of the
  vectors $ v(F) $ is bounded in terms of the size of $ \Delta $, i.e.\ that
  there are only finitely many choices for the direction vectors.
\end {proof}

\begin {proposition} \label {combtype-dim} \label {codim}
  For every combinatorial type $ \alpha $ occurring in $ \calM_n $ (resp.\ $
  \calM_{\Delta,n} $) the space $ \calM_n^\alpha $ (resp.\ $
  \calM_{\Delta,n}^\alpha $) is naturally an (unbounded) open convex polyhedron
  in a real vector space, i.e.\ a subset of a real vector space given by
  finitely many linear strict inequalities. Its dimension is as expected, i.e.\
  \begin {align*}
    \dim \calM_n^\alpha &= n-3-\codim \alpha \\
    \mbox {resp.} \quad
    \dim \calM_{\Delta,n}^\alpha &= |\Delta|-1+n-\codim \alpha.
  \end {align*}
\end {proposition}

\begin {proof}
  The first formula follows immediately from the combinatorial fact that a
  3-valent tropical curve with $n$ unbounded edges has exactly $ n-3 $ bounded
  edges: the space $ \calM_n^\alpha $ is simply parametrized by the
  lengths of all bounded edges, i.e.\ it is given as the subset of $
  \RR^{n-3-\codim \alpha} $ where all coordinates are positive.

  The statement about $ \calM_{\Delta,n}^\alpha $ follows in the same way,
  noting that a plane tropical curve in $ \calM_{\Delta,n} $ has $ |\Delta|+n
  $ unbounded edges and that we need two additional (unrestricted) parameters
  to describe translations, namely the coordinates of the image of a fixed
  ``root vertex'' $ V \in \Gamma^0 $.
\end {proof}

Ideally, one would of course like to make the spaces $ \calM_n $ and $
\calM_{\Delta,n} $ into tropical varieties themselves. Unfortunately, there is
however no general theory of tropical varieties yet. We will therefore work in
the category of polyhedral complexes, which will be sufficient for our
purposes.

\begin {definition}[Polyhedral complexes] \label {def-poly}
  Let $ X_1,\dots,X_N $ be (possibly unbounded) open convex polyhedra in real
  vector spaces. A \df {polyhedral complex} with cells $ X_1,\dots,X_N $ is a
  topological space $X$ together with continuous inclusion maps $ i_k:
  \overline {X_k} \to X $ such that $X$ is the disjoint union of the sets $
  i_k(X_k) $ and the ``coordinate changing maps'' $ i_k^{-1} \circ i_l $ are
  linear (where defined) for all $ k \neq l $. We will usually drop the
  inclusion maps $ i_k $ in the notation and say that the cells $ X_k $ are
  contained in $X$.

  The \df {dimension} $ \dim X $ of a polyhedral complex $X$ is the maximum
  of the dimensions of its cells. We say that $X$ is of \df {pure dimension}
  $ \dim X $ if every cell is contained in the closure of a cell of dimension
  $ \dim X $. A point of $X$ is said to be \df {in general position} if it is
  contained in a cell of dimension $ \dim X $.
\end {definition}

\begin {example} \label {ex-poly}
  The moduli spaces $ \calM_n $ and $ \calM_{\Delta,n} $ are polyhedral
  complexes of pure dimensions $ n-3 $ and $ |\Delta|-1+n $, respectively, with
  the cells corresponding to the combinatorial types. In fact, this follows
  from lemma \ref {fintypes} and proposition \ref {combtype-dim} together with
  the obvious remark that the boundaries of the cells $ \calM_n^\alpha $ (and $
  \calM_{\Delta,n}^\alpha $) can naturally be thought of as subsets of $
  \calM_n $ (resp.\ $ \calM_{\Delta,n} $) as well: they correspond to tropical
  curves where some of the bounded edges acquire zero length and finally
  vanish, leading to curves with vertices of higher valence. A tropical curve
  in $ \calM_n $ or $ \calM_{\Delta,n} $ is in general position if and only if
  it is 3-valent.
\end {example}


\section {Tropical multiplicities} \label {tropmult}

Having defined moduli spaces of abstract and plane tropical curves as
polyhedral complexes we will now go on and define morphisms between them.
Important properties of such morphisms will be their ``tropical''
multiplicities and degrees.

\begin {definition} \label {def-morph} \brk

  \vspace*{-3ex}

  \begin {enumerate}
  \item \label {def-morph-a}
    A \df {morphism} between two polyhedral complexes $X$ and $Y$ is a
    continuous map $ f:X \to Y $ such that for each cell $ X_i \subset X $ the
    image $ f(X_i) $ is contained in only one cell of $Y$, and $ f|_{X_i} $ is
    a linear map (of polyhedra).
  \item \label {def-morph-b}
    Let $ f: X \to Y $ be a morphism of polyhedral complexes of the same pure
    dimension, and let $ P \in X $ be a point such that both $P$ and $ f(P) $
    are in general position (in $X$ resp.\ $Y$). Then locally around $P$ the
    map $f$ is a linear map between vector spaces of the same dimension. We
    define the \df {multiplicity} $ \mult_f(P) $ of $f$ at $P$ to be the
    absolute value of the determinant of this linear map. Note that the
    multiplicity depends only on the cell of $X$ in which $P$ lies. We will
    therefore also call it the multiplicity of $f$ in this cell.
  \item \label {def-morph-c}
    Again let $ f: X \to Y $ be a morphism of polyhedral complexes of the same
    pure dimension. A point $ P \in Y $ is said to be \df {in $f$-general
    position} if $P$ is in general position in $Y$ and all points of $
    f^{-1}(P) $ are in general position in $X$. Note that the set of points in
    $f$-general position in $Y$ is the complement of a subset of $Y$ of
    dimension at most $ \dim Y-1 $; in particular it is a dense open subset.
    Now if $ P \in Y $ is a point in $f$-general position we define the \df
    {degree} of $f$ at $P$ to be
      \[ \deg_f(P) := \sum_{Q \in f^{-1}(P)} \mult_f(Q). \]
    Note that this sum is indeed finite: first of all there are only finitely
    many cells in $X$. Moreover, in each cell (of maximal dimension) of $X$
    where $f$ is not injective (i.e.\ where there might be infinitely many
    inverse image points of $P$) the determinant of $f$ is zero and hence so is
    the multiplicity for all points in this cell.

    Moreover, since $X$ and $Y$ are of the same pure dimension, the cones of
    $X$ on which $f$ is not injective are mapped to a locus of codimension at
    least 1 in $Y$. Thus the set of points in $f$-general position away from
    this locus is also a dense open subset of $Y$, and for all points in this
    locus we have that not only the sum above but indeed the fiber of $P$ is
    finite.
  \end {enumerate}
\end {definition}

\begin {remark} \label {coord-natural}
  Note that the definition of multiplicity in definition \ref {def-morph} \ref
  {def-morph-b} depends on the choice of coordinates on the cells of $X$ and
  $Y$. For the spaces $ \calM_n $ and $ \calM_{\Delta,n} $ (with cells $
  \calM_n^\alpha $ and $ \calM_{\Delta,n}^\alpha $) there were several equally
  natural choices of coordinates in the proof of proposition \ref
  {combtype-dim}: for graphs of a fixed combinatorial type we had to pick an
  ordering of the bounded edges and a root vertex. We claim that the
  coordinates for two different choices will simply differ by a linear
  isomorphism with determinant $ \pm 1 $. In fact, this is obvious for a
  relabeling of the bounded edges. As for a change of root vertex simply note
  that the difference $ h(V_2)-h(V_1) $ of the images of two vertices is given
  by $ \sum_F l([F]) \cdot v(F) $, where the sum is taken over the (unique)
  chain of flags leading from $ V_1 $ to $ V_2 $. This is obviously a linear
  combination of the lengths of the bounded edges, i.e.\ of the other
  coordinates in the cell. As these length coordinates themselves remain
  unchanged it is clear that the determinant of this change of coordinates is
  1. We conclude that \emph {the multiplicities and degrees of a morphism of
  polyhedral complexes whose source and/or target is a moduli space of abstract
  or plane tropical curves do not depend on any choices (of a root vertex or a
  labeling of the bounded edges)}.
\end {remark}

\begin {example} \label {ex-ev}
  For $ i \in \{1,\dots,n\} $ the \df {evaluation maps}
    \[ \ev_i: \calM_{\Delta,n} \to \RR^2, \quad
         (\Gamma,x_1,\dots,x_n,h) \mapsto h(x_i) \]
  are morphisms of polyhedral complexes. We denote the two coordinate functions
  of $ \ev_i $ by $ \ev_i^1, \ev_i^2 : \calM_{\Delta,n} \to \RR $ and the total
  evaluation map by $ \ev = \ev_1 \times \cdots \times \ev_n: \calM_{\Delta,n}
  \to \RR^{2n} $. Of course these maps are morphisms of polyhedral complexes as
  well.

  As a concrete example consider plane tropical curves of the following
  combinatorial types:
  \begin {enumerate} \parindent 0mm
  \item \label {ex-ev-a}
    For the combinatorial type

    \begin {center} \input {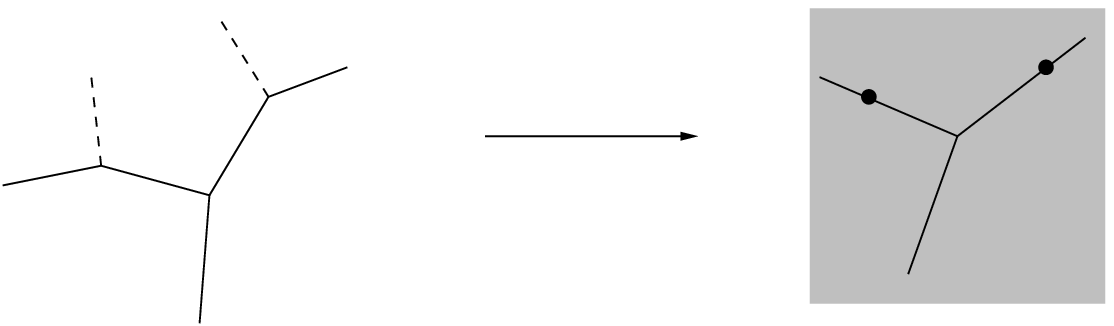} \end {center}

    we choose $V$ as the root vertex, say its image has coordinates $ h(V) =
    (a,b) $. There are two bounded edges with lengths $ l_i $ and direction
    vectors $ v_i = (v_{i,1},v_{i,2}) $ (counted from the root vertex) for $
    i=1,2 $. Then $ a,b,l_1,l_2 $ are the coordinates of $
    \calM_{\Delta,2}^\alpha $, and the evaluation maps are given by $ h(x_i) =
    h(V)+ l_i \cdot v_i = (a+l_i v_{i,1}, b+l_i v_{i,2} ) $. In particular,
    the total evaluation map $ \ev = \ev_1 \times \ev_2 $ is linear, and in
    the coordinates above its matrix is
      \[ \left( \begin {array}{cccc}
           1 & 0 & v_{1,1} &  0  \\
           0 & 1 & v_{1,2} &  0  \\
           1 & 0 &  0  & v_{2,1} \\
           0 & 1 &  0  & v_{2,2}
         \end {array} \right). \]
    An easy computation shows that the absolute value of the determinant of
    this matrix is $ \mult_{\ev} (\alpha) = | \det (v_1,v_2) | $. This is in
    fact the definition of the multiplicity $ \mult(V) $ of the vertex $V$ in
    \cite {Mi03} definition 4.15.
  \item \label {ex-ev-b}
    For the combinatorial type

    \begin {center} \input {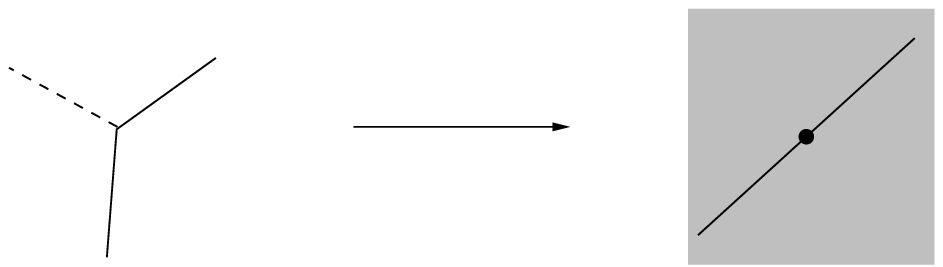} \end {center}

    the computation is even simpler: with the same reasoning as above the
    matrix of the evaluation map is just the $ 2 \times 2 $ unit matrix, and
    thus we get $ \mult_{\ev} (\alpha) = 1 $.
  \end {enumerate}
  Note that the entries of the matrices of evaluation maps will always be
  integers since the direction vectors of plane tropical curves lie in $ \ZZ^2
  $ by definition. In particular, multiplicities and degrees of evaluation maps
  will always be non-negative integers.
\end {example}

\begin {example} \label {ex-eval}
  Let $ n = |\Delta|-1 $, and consider the evaluation map $ \ev:
  \calM_{\Delta,n} \to \RR^{2n} $. Since both source and target of this map
  have dimension $ 2n $ we can consider the numbers
    \[ N_\Delta (\calP) := \deg_{\ev} (\calP) \in \ZZ_{\ge 0} \]
  for all points $ \calP \in \RR^{2n} $ in $ \ev $-general position.
  Note that these numbers are obviously just counting the tropical curves of
  degree $ \Delta $ through the points $ \calP $, where each such curve $C$ is
  counted with a certain multiplicity $ \mult_{\ev}(C) $. In the remaining part
  of this section we want to show how this multiplicity can be computed easily
  and that it is in fact the same as in definitions 4.15 and 4.16 of
  \cite {Mi03}.
\end {example}

\begin {definition} \label {def-rigid}
  Let $ C=(\Gamma,x_1,\dots,x_n,h) \in \calM_{\Delta,n} $ be a 3-valent plane
  tropical curve.
  \begin {enumerate}
  \item \label {def-rigid-a}
    A \df {string} of $C$ is a subgraph of $ \Gamma $ homeomorphic to $ \RR $
    (i.e.\ a ``path in $ \Gamma $ with two unbounded ends'') that does not
    intersect the closures $ \overline {x_i} $ of the marked points.
  \item \label {def-rigid-b}
    We say that (the combinatorial type of) $C$ is \df {rigid} if $ \Gamma $
    has no strings.
  \item \label {def-rigid-c}
    The \df {multiplicity} $ \mult (V) $ of a vertex $V$ of $C$ is defined
    to be $ |\det (v_1,v_2)| $, where $ v_1 $ and $ v_2 $ are two of the three
    direction vectors around $V$ (by the balancing condition it does not matter
    which ones we take here). The \df {multiplicity} $ \mult (C) $ of $C$ is
    the product of the multiplicities of all its vertices that are not adjacent
    to any marked point.
  \end {enumerate}
\end {definition}

\begin {remark} \label {rem-string}
  If $ C=(\Gamma,x_1,\dots,x_n,h) $ is a plane curve that contains a string $
  \Gamma' \subset \Gamma $ then there is a 1-parameter deformation of $C$ that
  moves the position of the string in $ \RR^2 $, but changes neither the images
  of the marked points nor the lines in $ \RR^2 $ on which the edges of $
  \Gamma \backslash \Gamma' $ lie. The following picture shows an example of
  (the image of) a plane 4-marked tropical curve with exactly one string $
  \Gamma' $ together with its corresponding deformation:

  \begin {center} \input {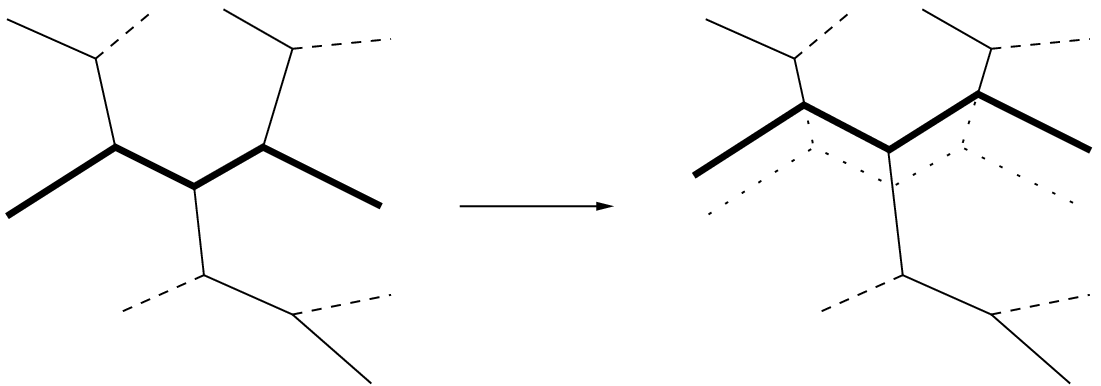} \end {center}
\end {remark}

\begin {remark} \label {rem-rigid}
  If $ C=(\Gamma,x_1,\dots,x_n,h) $ is an $n$-marked plane tropical curve of
  degree $ \Delta $ then the connected subgraph $ \Gamma \backslash \bigcup_i
  x_i $ has exactly $ |\Delta| $ unbounded ends. So if $ n < |\Delta|-1 $ there
  must be at least two unbounded ends that are still connected in $ \Gamma
  \backslash \bigcup_i \overline {x_i} $, i.e.\ there must be a string in $C$.
  If $ n=|\Delta|-1 $ then $C$ is rigid if and only if every connected
  component of $ \Gamma \backslash \bigcup_i \overline {x_i} $ has \emph
  {exactly} one unbounded end.
\end {remark}

\begin {proposition} \label {mult}
  Let $ n = |\Delta|-1 $. For any $n$-marked 3-valent plane tropical curve $C$
  we have
    \[ \mult_{\ev}(C) = \begin {cases}
         \mult (C) & \mbox {if $C$ is rigid,} \\
         0 & \mbox {otherwise,}
       \end {cases} \]
  where $ \mult(C) $ is as in definition \ref {def-rigid} \ref {def-rigid-c}.
\end {proposition}

\begin {proof}
  If $C$ is not rigid then by remark \ref {rem-string} it can be deformed with
  the images of the marked points fixed in $ \RR^2 $. This means that the
  evaluation map cannot be a local isomorphism and thus $ \mult_{\ev}(C)=0 $.
  We will therefore assume from now on that $C$ is rigid.

  We prove the statement by induction on the number $ k=2n-2 $ of bounded edges
  of $C$. The first cases $ k=0 $ and $ k=2 $ have been considered in example
  \ref {ex-ev}. So we can assume that $ k \ge 4 $. Choose any bounded edge $E$
  so that there is at least one bounded edge of $C$ to both sides of $E$. We
  cut $C$ along this edge into two halves $ C_1 $ and $ C_2 $. By extending
  the cut edge to infinity on both sides we can make $ C_1 $ and $ C_2 $ into
  plane tropical curves themselves:

  \begin {center} \input {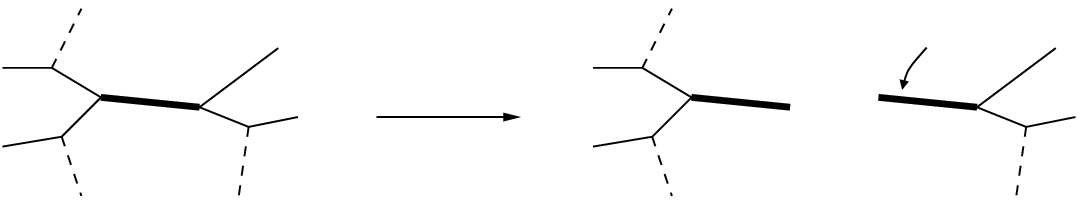} \end {center}

  (note that in this picture we have not drawn the map $h$ to $ \RR^2 $ but
  only the underlying abstract tropical curves). For $ i \in \{1,2\} $ we
  denote by $ n_i $ and $ k_i $ the number of marked points and bounded edges
  of $ C_i $, respectively. Of course we have $ n_1+n_2 = n $ and $ k_1+k_2 =
  k-1 = 2n-3 $.

  Assume first that $ k_1 \le 2n_1-3 $. As $ C_1 $ is 3-valent the total number
  of unbounded edges of $ C_1 $ is $ k_1 + 3 \le 2n_1 $; the number of unmarked
  unbounded edges is therefore at most $ n_1 $. This means that there must be
  at least one bounded connected component when we remove the closures of the
  marked points from $ C_1 $. The same is then true for $C$, i.e.\ by remark
  \ref {rem-rigid} $C$ is not rigid in contradiction to our assumption. By
  symmetry the same is of course true if $ k_2 \le 2n_2-3 $.

  The only possibility left is therefore $ k_1 = 2n_1-2 $ and $ k_2 = 2n_2-1 $
  (or vice versa). If we pick a root vertex in $ C_1 $ then in the matrix
  representation of the evaluation map we have $ 2n_1 $ coordinates in $
  \RR^{2n} $ (namely the images of the marked points on $ C_1 $) that depend
  on only $ 2+k_1 = 2n_1 $ coordinates (namely the root vertex and the lengths
  of the $ k_1 $ bounded edges in $ C_1 $). Hence the matrix has the form
    \[ \left( \begin {array}{c|c}
         A_1 & 0 \\ \hline
         * & A_2
       \end {array} \right) \]
  where $ A_1 $ and $ A_2 $ are square matrices of size $ 2n_1 $ and $ 2n_2 $,
  respectively. Note that $ A_1 $ is precisely the matrix of the evaluation
  map for $ C_1 $. As for $ A_2 $ its columns correspond to the lengths of $E$
  and the $ k_2 $ bounded edges of $ C_2 $, and its rows to the image points of
  the $ n_2 $ marked points on $ C_2 $. So if we consider the plane curve $
  \tilde C_2 $ obtained from $ C_2 $ by adding a marked point at a point $P$ on
  $E$ (see the picture above) and pick the vertex $P$ as the root vertex then
  the matrix for the evaluation map of $ \tilde C_2 $ is of the form
    \[ \left( \begin {array}{c|c}
         I_2 & 0 \\ \hline
         * & A_2
       \end {array} \right) \]
  where $ I_2 $ denotes the $ 2 \times 2 $ unit matrix and the two additional
  rows and columns correspond to the position of the root vertex. In particular
  this matrix has the same determinant as $ A_2 $. So we conclude that
    \[ \mult_{\ev}(C) = | \det A_1 \cdot \det A_2 | = \mult_{\ev_1} (C_1)
         \cdot \mult_{\ev_2}(\tilde C_2), \]
  where $ \ev_1 $ and $ \ev_2 $ denote the evaluation maps on $ C_1 $ and $
  \tilde C_2 $, respectively. The proposition now follows by induction, noting
  that $ C_1 $ and $ C_2 $ are rigid if $C$ is.
\end{proof}

\begin {remark} \label {rem-eval}
  By proposition \ref {mult} our numbers $ N_\Delta(\calP) $ are the
  same as the ones in \cite {Mi03}, and thus by the Correspondence Theorem
  (theorem 1 in \cite {Mi03}) the same as the corresponding complex numbers of
  stable maps. In particular they do not depend on $ \calP $ (as long as the
  points are in general position), and it is clear that the numbers $ N_d :=
  N_d(\calP) $ must satisfy Kontsevich's formula stated in the introduction. It
  is the goal of the rest of the paper to give an entirely tropical proof of
  this statement.
\end {remark}


\section {The forgetful maps} \label {forget}

We will now introduce the forgetful maps that have already been mentioned in
the introduction. As for the complex moduli spaces of stable maps there are
many such maps: given an $n$-marked plane tropical curve we can forget the map
to $ \RR^2 $, or some of the marked points, or both.

\begin {definition}[Forgetful maps] \label {def-forget}
  Let $ n \ge m $ be integers, and let $ C=(\Gamma,x_1,\dots,x_n,h) \in
  \calM_{\Delta,n} $ be an $n$-marked plane tropical curve.
  \begin {enumerate}
  \item \label {def-forget-a} (Forgetting the map and some points)
    Let $ C(m) $ be the minimal connected subgraph of $ \Gamma $ that
    contains the unbounded edges $ x_1,\dots,x_m $. Note that $ C(m) $
    cannot contain vertices of valence 1. So if we ``straighten'' the graph $
    C(m) $ at all 2-valent vertices (i.e.\ we replace the two adjacent edges
    and the vertex by one edge whose length is the sum of the lengths of the
    original edges) then we obtain an element of $ \calM_m $ that we denote by
    $ \ft_m (C) $.
  \item \label {def-forget-b} (Forgetting some points only)
    Let $ \tilde C(m) $ be the minimal connected subgraph of $ \Gamma $ that
    contains all unmarked ends as well as the marked points $ x_1,\dots,x_m $.
    Again $ \tilde C(m) $ cannot have vertices of valence 1. If we straighten $
    \tilde C(m) $ as in \ref {def-forget-a} we obtain an abstract tropical
    curve $ \tilde \Gamma $ with $ |\Delta|+m $ markings. Note that the
    restriction of $h$ to $ \tilde \Gamma $ still satisfies the requirements
    for a plane tropical curve, i.e.\ $ (\tilde \Gamma,x_1,\dots,x_m,h|_{\tilde
    \Gamma}) $ is an element of $ \calM_{\Delta,m} $. We denote it by $ \tilde
    {\ft}_m (C) $.
  \end {enumerate}
  It is obvious that the maps $ \ft_m: \calM_{\Delta,n} \to \calM_m $ and $
  \tilde {\ft}_m : \calM_{\Delta,n} \to \calM_{\Delta,m} $ defined in this way
  are morphisms of polyhedral complexes. We call them the \df {forgetful maps}
  (that keep only the first $m$ marked points resp.\ the first $m$ marked
  points and the map). Of course there are variations of the above maps: we can
  forget a given subset of the $n$ marked points that are not necessarily the
  last ones, or we can forget some points of an abstract tropical curve to
  obtain maps $ \calM_n \to \calM_m $.
\end {definition}

\begin {example} \label {ex-forget}
  For the plane tropical curve $C$ of example \ref {ex-stable-map} the graph $
  C(4) $ is simply the subgraph drawn in bold, and $ \ft_4(C) $ is the
  ``straightened version'' of this graph, i.e.\ the 4-marked tropical curve of
  type (A) in example \ref {ex-m4} with length parameter $l$ as indicated in
  the picture. Of course this length parameter is then also the local
  coordinate of $ \calM_4 $ if we want to represent the morphism $ \ft_4 $ of
  polyhedral complexes by a matrix, i.e.\ the matrix describing $ \ft_4 $ is
  the matrix with one row that has a 1 precisely at the column corresponding to
  the bounded edge marked $l$ (and zeroes otherwise).
\end {example}

The map that we need to consider for Kontsevich's formula is the following:

\begin {definition} \label {def-pi}
  Fix $ d \ge 2 $, and let $ n=3d $. We set
    \[ \pi := \ev_1^1 \times \ev_2^2 \times \ev_3 \times \cdots \times \ev_n
       \times \ft_4 : \calM_{d,n} \to \RR^{2n-2} \times \calM_4, \]
  i.e.\ $ \pi $ describes the first coordinate of the first marked point, the
  second coordinate of the second marked point, both coordinates of the other
  marked points, and the point in $ \calM_4 $ defined by the first four marked
  points. Obviously, $ \pi $ is a morphism of polyhedral complexes of pure
  dimension $ 2n-1 $.
\end {definition}

The central result of this section is the following proposition showing that
the degrees $ \deg_\pi (\calP) $ of $ \pi $ do not depend on the chosen point $
\calP $. Ideally this should simply follow from $ \pi $ being a ``morphism of
tropical varieties'' (and not just a morphism of polyhedral complexes). As
there is no such theory yet however we have to prove the independence of $
\calP $ directly.

\begin {proposition} \label {np-indep}
  The degrees $ \deg_\pi (\calP) $ do not depend on $ \calP $ (as long as $
  \calP $ is in $ \pi $-general position).
\end {proposition}

\begin {proof}
  It is clear that the degree of $ \pi $ is \emph {locally} constant on the
  subset of $ \RR^{2n-2} \times \calM_4 $ of points in $ \pi $-general position
  since at any curve that counts for $ \deg_\pi(\calP) $ with a non-zero
  multiplicity the map $ \pi $ is a local isomorphism. Recall that the points
  in $ \pi $-general position are the complement of a polyhedral complex of
  codimension 1, i.e.\ they form a finite number of top-dimensional regions
  separated by ``walls'' that are polyhedra of codimension 1. Hence to show
  that $ \deg_\pi $ is globally constant it suffices to consider a general
  point on such a wall and to show that $ \deg_\pi $ is locally constant at
  these points too. Such a general point on a wall is simply the image under $
  \pi $ of a general plane tropical curve $C$ of a combinatorial type of
  codimension 1. So we simply have to check that $ \deg_\pi $ is locally
  constant around such a point $ C \in \calM_{\Delta,n} $.

  By definition a combinatorial type $ \alpha $ of codimension 1 has exactly
  one 4-valent vertex $V$, with all other vertices being 3-valent. Let $
  E_1,\dots,E_4 $ denote the four (bounded or unbounded) edges around $V$.
  There are precisely 3 combinatorial types $ \alpha_1,\alpha_2,\alpha_3 $ that
  have $ \alpha $ in their boundary, as indicated in the following local
  picture:

  \begin {center} \input {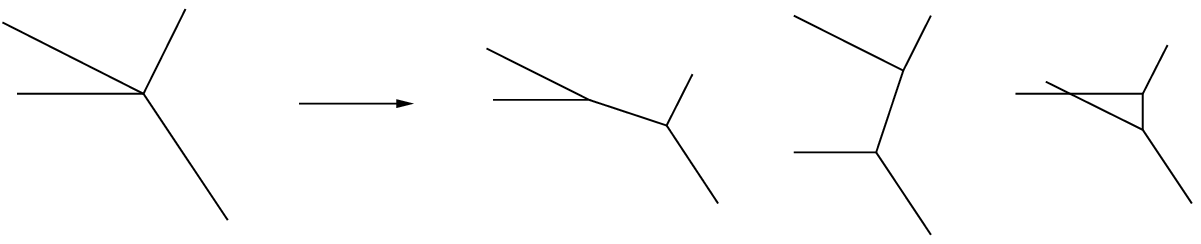} \end {center}

  Let us assume first that all four edges $ E_i $ are bounded. We denote their
  lengths by $ l_i $ and their directions (pointing away from $V$) by $ v_i $.
  To set up the matrices of $ \pi $ we choose the root vertex $V$ in $ \alpha_i
  $ as in the picture. We denote its image by $ w \in \RR^2 $.

  The following table shows the relevant parts of the matrices $ A_i $
  of $ \pi $ for the three combinatorial types $ \alpha_i $. Each matrix
  contains the first block of columns (corresponding to the image $w$ of the
  root vertex and the lengths $ l_i $ of the edges $ E_i $) and the $i$-th of
  the last three columns (corresponding to the length of the newly added
  bounded edge). The columns corresponding to the other bounded edges are not
  shown; it suffices to note here that they are the same for all three
  matrices. All rows but the last one correspond to the images in $ \RR^2 $ of
  the marked points; we get different types of rows depending on via which edge
  $ E_i $ this marked point can be reached from $V$. For the marked points $
  x_i $ with $ i \ge 3 $ we use both coordinates in $ \RR^2 $ (hence one row
  in the table below corresponds to two rows in the matrix), for $ x_1 $ only
  the first and for $ x_2 $ only the second coordinate. The last row
  corresponds to the coordinate in $ \calM_4 $ as in example \ref {ex-forget}.
  In the following table $ I_2 $ denotes the $ 2 \times 2 $ unit matrix, and
  each $*$ and $**$ stands for 0 or 1 (see below).
    \[ \begin {array}{l|ccccc|c|c|c}
         & w & l_1 & l_2 & l_3 & l_4 & l^{\alpha_1}
           & l^{\alpha_2} & l^{\alpha_3} \\ \hline
         \mbox {points behind $ E_1 $} &
           I_2 & v_1 & 0 & 0 & 0 & 0 & 0 & 0 \\
         \mbox {points behind $ E_2 $} &
           I_2 & 0 & v_2 & 0 & 0 & v_2+v_3 & 0 & v_2+v_4 \\
         \mbox {points behind $ E_3 $} &
           I_2 & 0 & 0 & v_3 & 0 & v_2+v_3 & v_3+v_4 & 0 \\
         \mbox {points behind $ E_4 $} &
           I_2 & 0 & 0 & 0 & v_4 & 0 & v_3+v_4 & v_2+v_4 \\
         \mbox {coordinate of $ \calM_4 $} &
           0 & * & * & * & * & ** & ** & **
       \end {array} \]
  To look at these matrices (in particular at the entries marked $*$) further
  we will distinguish several cases depending on how many of the edges $
  E_1,\dots,E_4 $ of $C$ are contained in the subgraph $ C(4) $ of definition
  \ref {def-forget}:
  \begin {enumerate} \parindent 0mm \parskip 0.5ex
  \item \label {np-indep-a}
    4 edges: Then $ \ft_4(C) $ is the curve (D) of example \ref {ex-m4}, and
    the three types $ \alpha_1,\alpha_2,\alpha_3 $ are mapped precisely to the
    three other types (A), (B), (C) of $ \calM_4 $ by $ \ft_4 $, i.e.\ to the
    three cells of $ \RR^{2n-2} \times \calM_4 $ around the wall by $ \pi $.
    For these three types the length parameter in $ \calM_4 $ is simply the one
    newly inserted edge. Hence the entries $*$ in the matrix above are all 0,
    whereas the entries $**$ are all 1. It follows that the three matrices $
    A_1,A_2,A_3 $ have a 1 as the bottom right entry and all zeroes in the
    remaining places of the last row. Their determinants therefore do not
    depend on the last column. But this is the only column that differs for the
    three matrices, i.e.\ $ A_1,A_2, $ and $ A_3 $ all have the same
    determinant. It follows by definition that $ \deg_\pi $ is locally constant
    around $C$. This completes the proof of the proposition in this case.
  \item \label {np-indep-b}
    3 edges: The following picture shows what the combinatorial types $
    \alpha $, $ \alpha_1 $, $ \alpha_2 $, $ \alpha_3 $ look like locally around
    the vertex $V$ in this case. As in example \ref {ex-stable-map} we have
    drawn the edges belonging to $ C(4) $ in bold.

    \begin {center} \input {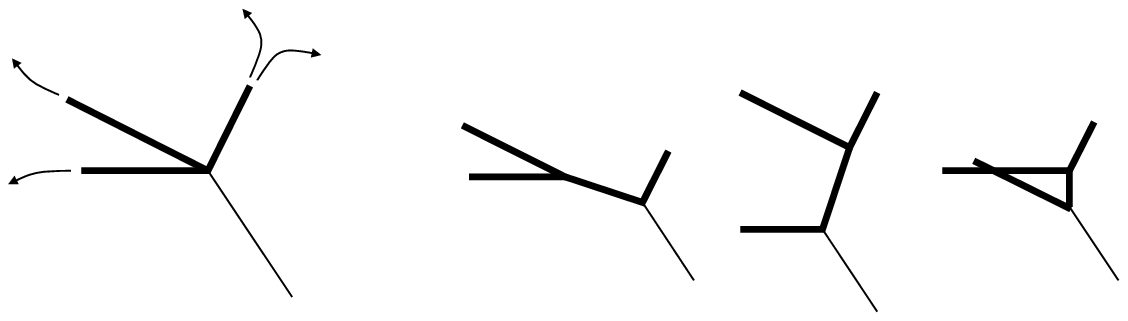} \end {center}

    We see that exactly one edge $ E_i $ (namely $ E_2 $ in the example above)
    counts towards the length parameter in $ \calM_4 $, and that the newly
    inserted edge counts towards this length parameter in exactly one of the
    combinatorial types $ \alpha_i $ (namely $ \alpha_1 $ in the example
    above). Hence in the table showing the matrices $ A_i $ above exactly one
    of the entries $*$ and exactly one of the entries $**$ is 1, whereas the
    others are 0.
  \item \label {np-indep-c}
    2 edges: There are two possibilities in this case. If $V$ is a point in $
    C(4) $ corresponding to an interior point of the bounded edge in $ \ft_4(C)
    $ then an analysis completely analogous to that in \ref {np-indep-b} shows
    that exactly 2 of the entries $*$ and also 2 of the entries $**$ above are
    1, whereas the others are zero. If on the other hand $V$ corresponds to an
    interior point of an unbounded edge in $ \ft_4(C) $ then all entries $*$
    and $**$ above are 0.
  \item \label {np-indep-d}
    fewer than 2 edges: As it is not possible that exactly one of the edges $
    E_i $ is contained in $ C(4) $ we must then have that there is no such
    edge, and consequently that all entries $*$ and $**$ above are 0.
  \end {enumerate}
  Summarizing, we see in all remaining cases \ref {np-indep-b}, \ref
  {np-indep-c}, and \ref {np-indep-d} that there are equally many entries
  $**$ equal to 1 as there are entries $*$ equal to 1. So using the linearity
  of the determinant in the column corresponding to the new bounded edge we
  get that $ \det A_1 + \det A_2 + \det A_3 $ is equal to the determinant of
  the matrix whose corresponding entries are
    \[ \begin {array}{l|ccccc|c}
         & w & l_1 & l_2 & l_3 & l_4 & l \\ \hline
         \mbox {points behind $ E_1 $} &
           I_2 & v_1 & 0 & 0 & 0 & 0 \\
         \mbox {points behind $ E_2 $} &
           I_2 & 0 & v_2 & 0 & 0 & 2v_2+v_3+v_4 \\
         \mbox {points behind $ E_3 $} &
           I_2 & 0 & 0 & v_3 & 0 & 2v_3+v_2+v_4 \\
         \mbox {points behind $ E_4 $} &
           I_2 & 0 & 0 & 0 & v_4 & 2v_4+v_2+v_3 \\
         \mbox {coordinate of $ \calM_4 $} &
           0 & * & * & * & * & **
       \end {array} \]
  where $**$ is now the sum of the four entries marked $*$. If we now subtract
  the four $ l_i $-columns and add $ v_1 $ times the $w$-columns from the last
  one then all entries in the last column vanish (note that $ v_1+v_2+v_3+v_4=0
  $ by the balancing condition). So we conclude that
    \[ \det A_1 + \det A_2 + \det A_3 = 0. \tag {$1$} \]
  For a given $ i \in \{1,2,3\} $ let us now determine whether the
  combinatorial type $ \alpha_i $ occurs in the inverse image of a fixed point
  $ \calP $ near the wall. We may assume without loss of generality that the
  multiplicity of $ \alpha_i $ is non-zero since other types are irrelevant for
  the statement of the proposition. So the restriction $ \pi_i $ of $ \pi $ to
  $ \calM_{\Delta,n}^{\alpha_i} $ is given by the invertible matrix $ A_i $.
  There is therefore at most one inverse image point in $ \pi_i^{-1}(\calP) $,
  which would have to be the point with coordinates $ A_i^{-1} \cdot \calP $.
  In fact, this point exists in $ \calM_{\Delta,n}^{\alpha_i} $ if and only if
  all coordinates of $ A_i^{-1} \cdot \calP $ corresponding to lengths of
  bounded edges are positive. By continuity this is obvious for all edges
  except the newly added one since in the boundary curve $C$ all these edges
  had positive length. We conclude that there is a point in $ \pi_i^{-1}(\calP)
  $ if and only if the last coordinate (corresponding to the length of the
  newly added edge) of $ A_i^{-1} \cdot \calP $ is positive. By Cramer's rule
  this last coordinate is $ \det \tilde A_i / \det A_i $, where $ \tilde A_i $
  denotes the matrix $ A_i $ with the last column replaced by $ \calP $. But
  note that $ \tilde A_i $ does not depend on $i$ since the last column was the
  only one where the $ A_i $ differ. Hence whether there is a point in $
  \pi_i^{-1}(\calP) $ or not depends solely on the sign of $ \det A_i $: either
  there are such inverse image points for exactly those $i$ where $ \det A_i $
  is positive, or exactly for those $i$ where $ \det A_i $ is negative. But by
  $(1)$ the sum of the absolute values of the determinants satisfying this
  condition is the same in both cases. This means that $ \deg_\pi $ is locally
  constant around $C$.

  Strictly speaking we have assumed in the above proof that all edges $ E_i $
  are bounded. It is very easy however to adapt these arguments to the other
  cases: if an edge $ E_i $ is not bounded then there is no coordinate $ l_i $
  corresponding to its length, but neither are there marked points that can be
  reached from $V$ via $ E_i $. We leave it as an exercise to check that the
  above proof still holds in this case with essentially no modifications.
\end {proof}


\section {Kontsevich's formula} \label {kontsevich}

We have just shown that the degrees of the morphism $ \pi: \calM_{d,n} \to
\RR^{2n-2} \times \calM_4 $ of definition \ref {def-pi} do not depend on the
point chosen in the target. We now want to apply this result by equating the
degrees for two different points, namely two points where the $ \calM_4
$-coordinate is very large, but corresponds to curves of type (A) or (B) in
example \ref {ex-m4}. We will first prove that a very large length in $ \calM_4
$ requires the curves to acquire a contracted bounded edge.

\begin {proposition} \label {m4-bounded}
  Let $ d \ge 2 $ and $ n=3d $, and let $ \calP \in \RR^{2n-2} \times \calM_4 $
  be a point in $ \pi $-general position whose $ \calM_4 $-coordinate is very
  large (i.e.\ it corresponds to a 4-marked curve of type (A), (B), or (C) as
  in example \ref {ex-m4} with a very large length $l$). Then every plane
  tropical curve $ C \in \pi^{-1}(\calP) $ with $ \mult_{\pi}(C) \neq 0 $ has
  a contracted bounded edge.
\end {proposition}

\begin {proof}
  We have to show that the set of all points $ \ft_4 (C) \in \calM_4 $ is
  bounded in $ \calM_4 $, where $C$ runs over all curves in $ \calM_{d,n} $
  with non-zero $ \pi $-multiplicity that have no contracted bounded edge and
  satisfy the given incidence conditions at the marked points. As there are
  only finitely many combinatorial types by lemma \ref {fintypes} we can
  restrict ourselves to curves of a fixed (but arbitrary) combinatorial type
  $ \alpha $. Since $ \calP $ is in $ \pi $-general position we can assume that 
  the codimension of $ \alpha $ is 0, i.e.\ that the curve is 3-valent.

  Let $ C' \in \calM_{d,n-2} $ be the curve obtained from $C$ by forgetting the
  first two marked points as in definition \ref {def-forget}. We claim that $
  C' $ has exactly one string (see definition \ref {def-rigid} \ref
  {def-rigid-a}). In fact, $ C' $ must have at least one string by remark \ref
  {rem-rigid} since $C'$ has less than $ 3d-1=n-1 $ marked points. On the other
  hand, if $C'$ had at least two strings then by remark \ref {rem-string} $C'$
  would move in an at least 2-dimensional family with the images of $
  x_3,\dots,x_n $ fixed. It follows that $C$ moves in an at least 2-dimensional
  family as well with the first coordinate of $ x_1 $, the second of $ x_2 $,
  and both of $ x_3,\dots,x_n $ fixed. As $ \calM_4 $ is one-dimensional this
  means that $C$ moves in an at least 1-dimensional family with the image point
  under $ \pi $ fixed. Hence $ \pi $ is not a local isomorphism, i.e.\ $
  \mult_\pi(C)=0 $ in contradiction to our assumptions.

  So let $ \Gamma' $ be the unique string in $C'$. The deformations of $C'$
  with the given incidence conditions fixed are then precisely the ones of the
  string described in remark \ref {rem-string}. Note that the edges adjacent to
  $ \Gamma' $ must be bounded since otherwise we would have two strings. So if
  there are edges adjacent to $ \Gamma' $ to both sides of $ \Gamma' $ as in
  picture (a) below (note that there are no contracted bounded edges by
  assumption) then the deformations of $C'$ with the combinatorial type and the
  incidence conditions fixed are bounded on both sides. For the deformations of
  $C$ with its combinatorial type and the incidence conditions fixed this means
  that the lengths of all inner edges are bounded except possibly the edges
  adjacent to $ x_1 $ and $ x_2 $. This is sufficient to ensure that the
  image of these curves under $ \ft_4 $ is bounded in $ \calM_4 $ as well.

  \begin {center} \input {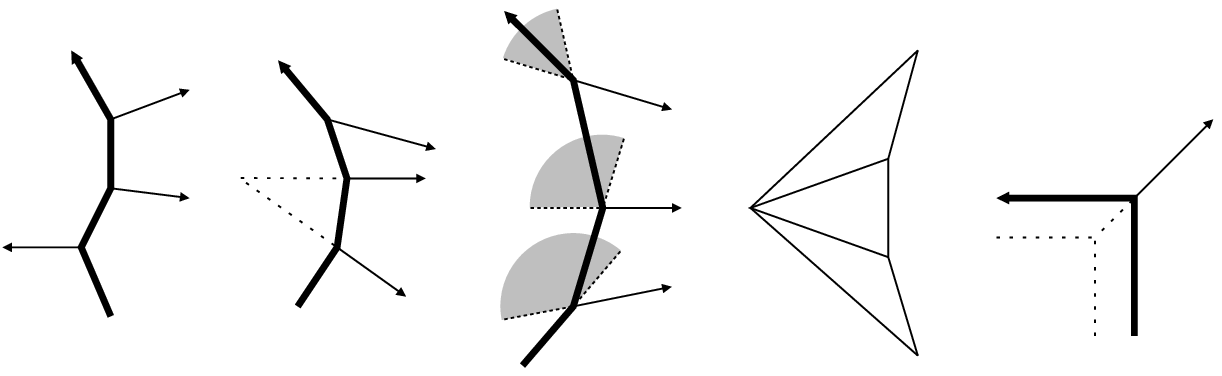} \end {center}

  Hence we are only left with the case when all adjacent edges of $ \Gamma' $
  are on the same side of $ \Gamma' $, say after picking an orientation of $
  \Gamma $ on the right side as in picture (b) above. Label the edges (resp.\
  their direction vectors) of $ \Gamma' $ by $ v_1,\dots,v_k $ and the adjacent
  edges of the curve by $ w_1,\dots,w_{k-1} $ as in the picture. As above the
  movement of $C'$ to the right within its combinatorial type is bounded. If
  one of the directions $ w_{i+1} $ is obtained from $ w_i $ by a left turn (as
  it is the case for $ i=1 $ in the picture) then the edges $ w_i $ and $
  w_{i+1} $ meet to the left of $ \Gamma' $. This restricts the movement of
  $C'$ to the left within its combinatorial type too since the corresponding
  edge $ v_{i+1} $ then shrinks to zero. We can then conclude as in case (a)
  above that the image of these curves under $ \ft_4 $ is bounded.

  We can therefore assume that for all $i$ the direction $ w_{i+1} $ is
  either the same as $ w_i $ or obtained from $ w_i $ by a right turn as in
  picture (c). The balancing condition then ensures that for all $i$ both the
  directions $ v_{i+1} $ and $ -w_{i+1} $ lie in the angle between $ v_i $ and
  $ -w_i $ (shaded in the picture above). It follows that all directions $ v_i
  $ and $ -w_i $ lie within the angle between $ v_1 $ and $ -w_1 $. In
  particular the string $ \Gamma' $ cannot have any self-intersections in $
  \RR^2 $. We can therefore pass to the (local) dual picture (d) (see e.g.\
  \cite {Mi03} section 3.4) where the edges dual to $ w_i $ correspond to a
  concave side of the polygon whose other two edges are the ones dual to $ v_1
  $ and $ v_k $. In other words, the intersection points of the edges dual to $
  w_{i-1} $ and $ w_i $ must be in the interior of the triangle spanned by the
  edges dual to $ v_1 $ and $ v_k $ for all $ 1 < i < k $.

  But note that both $ v_1 $ and $ v_k $ must be $ (-1,0) $, $ (0,-1) $, or $
  (1,1) $ since they are outer directions of a curve of degree $d$.
  Consequently, their dual edges have to be among the vectors $ \pm (1,0) $,
  $ \pm (0,1) $, $ \pm (1,-1) $. But any triangle spanned by two of these
  vectors has area (at most) $ \frac 12 $ and thus does not admit any integer
  interior points. It follows that intersection points of the dual edges of $
  w_{i-1} $ and $ w_i $ as above cannot exist and therefore that $ k=2 $,
  i.e.\ that the string consists just of the two unbounded ends $ v_1 $ and
  $ v_2 $ that are connected to the rest of the curve by exactly one internal
  edge $ w_1 $. It must therefore look as in picture (e).

  In this case the movement of the string is indeed not bounded to the left.
  Note that then $ w_1 $ is the only internal edge whose length is not
  bounded within the deformations of $C'$ since the rest of the curve (not
  shown in picture (e)) does not move at all. But we will show that this
  unbounded length of $ w_1 $ cannot count towards the length parameter in
  $ \calM_4 $ for the deformations of $C$: first of all this would require
  two of the marked points $ x_1,\dots,x_4 $ to lie on $ v_1 $ or $ v_2 $ for
  all curves in the deformation, but of course with $ v_1 $ and $ v_2 $
  forming a string we cannot have $ x_3 $ or $ x_4 $ (where we impose point
  conditions) on them. Hence we would have to have $ x_1 $ and $ x_2 $ (that we
  require to lie on a vertical line $ L_1 $ resp.\ a horizontal line $ L_2 $)
  somewhere on $ v_1 $ and $ v_2 $. But the following picture shows that for
  all three possibilities for $ v_1 $ and $ v_2 $ the union of the edges $ v_1
  $ and $ v_2 $ (drawn in bold) finally becomes disjoint from at least one of
  the lines $ L_1 $ and $ L_2 $ as the length of $ w_1 $ increases:

  \begin {center} \input {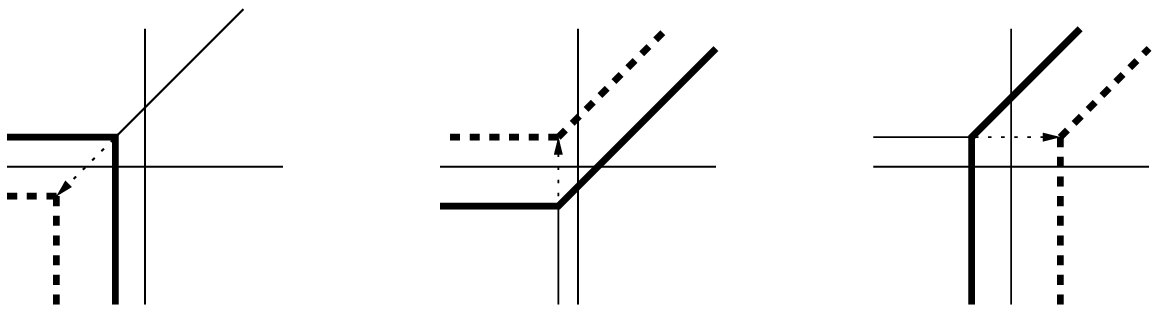} \end {center}

  This means that we cannot have both $ x_1 $ and $ x_2 $ on the union of $
  v_1 $ and $ v_2 $ as the length of $ w_1 $ increases. Consequently, we
  cannot get unbounded length parameters in $ \calM_4 $ in this case either.
  This finishes the proof of the proposition.
\end {proof}

\begin {remark} \label {rem-reducible}
  Let $ C=(\Gamma,x_1,\dots,x_n,h) $ be a plane tropical curve with a
  contracted bounded edge $E$, and assume that there is at least one more
  bounded edge to both sides of $E$. Then in the same way as in
  the proof of proposition \ref {mult} we can split $ \Gamma $ at $E$ into two
  graphs $ \Gamma_1 $ and $ \Gamma_2 $, making the edge $E$ into a contracted
  \emph {unbounded} edge on both sides. By restricting $h$ to these graphs we
  obtain two new plane tropical curves $ C_1 $ and $ C_2 $. The marked points $
  x_1,\dots,x_n $ obviously split up onto $ C_1 $ and $ C_2 $; in addition
  there is one more marked point $P$ resp.\ $Q$ on both curves that corresponds
  to the newly added contracted unbounded edge. If $C$ is a curve of degree $d$
  then (by the balancing condition) $ C_1 $ and $ C_2 $ are of some degrees $
  d_1 $ and $ d_2 $ with $ d_1+d_2=d $.

  \begin {center} \input {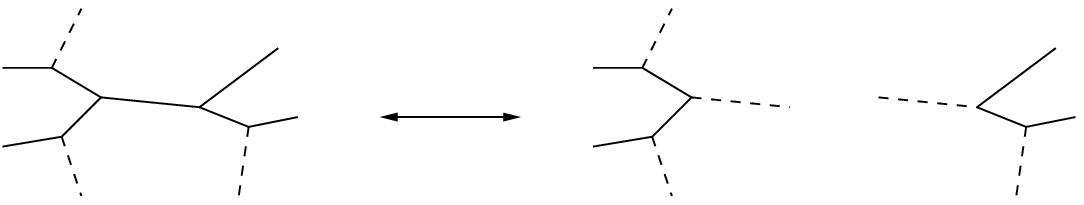} \end {center}

  We will say in this situation that $C$ is obtained by \df {glueing} $ C_1 $
  and $ C_2 $ along the identification $ P=Q $, and that $C$ is a \df
  {reducible} plane tropical curve that can be decomposed into $ C_1 $ and $
  C_2 $. For the image we obviously have $ h(\Gamma) = h(\Gamma_1) \cup
  h(\Gamma_2) $, so when considering embedded plane tropical curves $C$ is in
  fact just the union of the two curves $ C_1 $ and $ C_2 $ of smaller degree
  (see example \ref {ex-conic-red}).
\end {remark}

\begin {lemma} \label {set-reducible}
  Let $ \calP =(a,b,p_3,\dots,p_n,z) \in \RR^{2n-2} \times \calM_4 $ be a point
  in $ \pi $-general position such that $ z \in \calM_4 $ is of type (A) (see
  example \ref {ex-m4}) with a very large length parameter. Then for every
  plane tropical curve $C$ in $ \pi^{-1}(\calP) $ with non-zero $ \pi
  $-multiplicity we have exactly one of the following cases:
  \begin {enumerate}
  \item \label {set-reducible-a}
    $ x_1 $ and $ x_2 $ are adjacent to the same vertex (that maps to $ (a,b) $
    under $h$);
  \item \label {set-reducible-b}
    $C$ is reducible and decomposes uniquely into two components $ C_1 $ and $
    C_2 $ of some degrees $ d_1 $ and $ d_2 $ with $ d_1+d_2=d $ such that the
    marked points $ x_1 $ and $ x_2 $ are on $ C_1 $, the points $ x_3 $ and $
    x_4 $ are on $ C_2 $, and exactly $ 3d_1-1 $ of the other points $ x_5,
    \dots,x_n $ are on $ C_1 $.
  \end {enumerate}
\end {lemma}

\begin {proof}
  By proposition \ref {m4-bounded} any curve $ C \in \pi^{-1}(\calP) $ with
  non-zero $ \pi $-multiplicity has at least one contracted bounded edge. In
  fact, $C$ must have \emph {exactly} one such edge: if $C$ had at least 2
  contracted bounded edges then there would be $ 2n-2 $ coordinates in the
  target of $ \pi $ (namely the evaluation maps) that depend on only $ 2n-3 $
  variables (namely the root vertex and the lengths of all but 2 of the $ 2n-3
  $ bounded edges), hence we would have $ \mult_\pi(C) = 0 $.

  So let $E$ be the unique contracted bounded edge of $C$. Note that $E$ must
  be contained in the subgraph $ C(4) $ of definition \ref {def-forget} \ref
  {def-forget-a} since otherwise we could not have a very large length
  parameter in $ \calM_4 $. As the point $z$ is of type (A) we conclude that $
  x_1 $ and $ x_2 $ must be to one side, and $ x_3 $ and $ x_4 $ to the other
  side of $E$. Denote these sides by $ C_1 $ and $ C_2 $, respectively.

  If there are no bounded edges in $ C_1 $ then $C$ is not reducible as in
  remark \ref {rem-reducible}. Instead $ C_1 $ consists only of $E$, $ x_1
  $, and $ x_2 $, i.e.\ we are then in case \ref {set-reducible-a}. The
  evaluation conditions then require that all of $ C_1 $ must be mapped to the
  point $ (a,b) $. Note that it is not possible that there are no bounded edges
  in $ C_2 $ since this would require $ x_3 $ and $ x_4 $ to map to the same
  point in $ \RR^2 $.

  We are left with the case when there are bounded edges to both sides of $E$.
  In this case $C$ is reducible as in remark \ref {rem-reducible}, so we are
  in case \ref {set-reducible-b}. In this case $ x_1 $ and $ x_2 $ cannot be
  adjacent to the same vertex since this would require another contracted edge
  by the balancing condition. Now let $ n_1 $ and $ n_2 $ be the number of
  marked points $ x_5,\dots,x_n $ on $ C_1 $ resp.\ $ C_2 $; we have to show
  that $ n_1 = 3d_1-1 $ and $ n_2 = 3d_2-3 $. So assume that $ n_1 \ge 3d_1 $.
  Then at least $ 2n_1+2 \ge 3d_1 +n_1+2 $ of the coordinates of $ \pi $ (the
  images of the $ n_1 $ marked points as well as the first image coordinate of
  $ x_1 $ and the second of $ x_2 $) would depend on only $ 3d_1+n_1+1 $
  coordinates (2 for the root vertex and one for each of the $ 3d_1+(n_1+2)-3 $
  bounded edges), leading to a zero $ \pi $-multiplicity. Hence we conclude
  that $ n_1 \le 3d_1-1 $. The same argument shows that $ n_2 \le 3d_2-3 $, so
  as the total number of points is $ n_1+n_2 = n-4 = (3d_1-1)+(3d_2-3) $ it
  follows that we must have equality.
\end {proof}

\begin {remark} \label {converse-reducible}
  In fact, the following ``converse'' of lemma \ref {set-reducible} is also
  true: as above let $ \calP =(a,b,p_3,\dots,p_n,z) \in \RR^{2n-2} \times
  \calM_4 $ be a point in $ \pi $-general position such that $ z \in \calM_4 $
  is of type (A) (see example \ref {ex-m4}) with a very large length parameter.
  Now let $ C_1 $ and $ C_2 $ be two (unmarked) plane tropical curves of
  degrees $ d_1 $ and $ d_2 $ with $ d_1+d_2=d $ such that the image of $ C_1 $
  passes through $ L_1:= \{ (x,y);\; x=a\} $, $ L_2 := \{ (x,y);\; y=b \} $,
  and $ 3d_1-1 $ of the points $ p_5,\dots,p_n $, and the image of $ C_2 $
  through $ p_3 $, $ p_4 $, and the other $ 3d_2-3 $ of the points $
  p_5,\dots,p_n $.

  Then for each choice of points $ P \in C_1 $ and $ Q \in C_2 $ that map to
  the same image point in $ \RR^2 $, and for each choice of points $
  x_1,\dots,x_n $ on $ C_1 $ and $ C_2 $ that map to $ L_1 $, $ L_2 $, $ p_3,
  \dots,p_n $, respectively, we can make $ C_1 $ and $ C_2 $ into marked plane
  tropical curves and glue them together to a single reducible $n$-marked curve
  $C$ in $ \pi^{-1}(\calP) $ as in remark \ref {rem-reducible} (the length of
  the one contracted edge is determined by $z$).

  As $ \calP $ was assumed to be in $ \pi $-general position we can never
  construct a curve $C$ in this way that is not 3-valent. In particular this
  means for example that $ C_1 $ and $ C_2 $ are guaranteed to be 3-valent
  themselves. Moreover, a point that is in the image of both $ C_1 $ and $ C_2
  $ cannot be a vertex of either curve. In particular, it is not possible that
  $ C_1 $ and $ C_2 $ share a common line segment in $ \RR^2 $. In the same way
  we see that the image of $ C_1 $ cannot meet $ L_1 $ or $ L_2 $ in a vertex
  or have a line segment in common with $ L_1 $ or $ L_2 $, and cannot meet $
  L_1 \cap L_2 $ at all.

  Summarizing, we see that after choosing the two curves $ C_1 $ and $ C_2 $ as
  well as the points $ x_1,\dots,x_n,P,Q $ on them there is a unique curve in $
  \pi^{-1}(\calP) $ obtained from this data. So if we want to compute the
  degree of $ \pi $ and have to sum over all points in $ \pi^{-1}(\calP) $ then
  for the curves of type \ref {set-reducible-b} in lemma \ref {set-reducible}
  we can as well sum over all choices of $ C_1 $, $ C_2 $, $ x_1,\dots,x_n,P,Q
  $ as above.
\end {remark}

Before we can actually do the summation we still have to compute the
multiplicity of $ \pi $ at the curves in $ \pi^{-1}(\calP) $:

\begin {proposition} \label {comp-mult}
  With notations as in lemma \ref {set-reducible} and remark \ref
  {converse-reducible} let $C$ be a point in $ \pi^{-1}(\calP) $. Then
  \begin {enumerate}
  \item \label {comp-mult-a}
    if $C$ is of type \ref {set-reducible-a} as in lemma \ref {set-reducible}
    its $ \pi $-multiplicity is $ \mult_{\ev}(C') $, where $C'$ denotes the
    curve obtained from $C$ by forgetting $ x_1 $, and $ \ev $ is the
    evaluation at the $ 3d-1 $ points $ x_2,\dots,x_n $;
  \item \label {comp-mult-b}
    if $C$ is of type \ref {set-reducible-b} as in lemma \ref {set-reducible}
    its $ \pi $-multiplicity is
    \[ \mult_\pi (C) = \mult_{\ev}(C_1) \cdot \mult_{\ev}(C_2) \cdot
         (C_1 \cdot C_2)_{P=Q} \cdot (C_1 \cdot L_1)_{x_1} \cdot (C_1 \cdot
         L_2)_{x_2}, \]
    where $ \mult_{\ev}(C_i) $ denotes the multiplicities of the evaluation map
    at the $ 3d_i-1 $ points of $ x_3,\dots,x_n $ that lie on the respective
    curve, and $ (C' \cdot C'')_P $ denotes the intersection multiplicity of
    the tropical curves $ C' $ and $ C'' $ at the point $ P \in C' \cap C'' $
    (see \cite {RST03} section 4), i.e.\ $ |\det (v',v'')| $ where $ v' $ and $
    v'' $ are the direction vectors of $C'$ and $C''$ at $P$. In particular,
    $ (C_1 \cdot L_i)_{x_i} $ is simply the first resp.\ second coordinate of
    the direction vector of $ C_1 $ at $ x_i $ for $ i \in \{1,2\} $.
  \end {enumerate}
\end {proposition}

\begin {proof}
  We simply have to set up the matrix for $ \pi $ and compute its determinant.
  First of all note that in both cases \ref {comp-mult-a} and \ref
  {comp-mult-b} the length of the contracted bounded edge is irrelevant for all
  evaluation maps and contributes with a factor of 1 to the $
  \calM_4 $-coordinate of $ \pi $. Hence the column of $ \pi $ corresponding to
  the contracted bounded edge has only one entry 1 and all others zero. To
  compute its determinant we may therefore drop both the $ \calM_4 $-row and
  the column corresponding to the contracted bounded edge.

  In case \ref {comp-mult-a} the matrix obtained this way is then exactly the
  same as if we had only one marked point instead of $ x_1 $ and $ x_2 $ and
  evaluate this point for both coordinates in $ \RR^2 $ (instead of evaluating
  $ x_1 $ for the first and $ x_2 $ for the second). This proves \ref
  {comp-mult-a}.

  For \ref {comp-mult-b} let us first consider the marked point $ x_1 $ where
  we only evaluate the first coordinate. Let $ E_1 $ and $ E_2 $ be the two
  adjacent edges and assume first that both of them are bounded. Denote their
  common direction vector by $ v=(v^1,v^2) $ and their lengths by $ l_1, l_2 $.
  Assume that the root vertex is on the $ E_1 $-side of $ x_1 $. Then the
  entries of the matrix for $ \pi $ corresponding to $ l_1 $ and $ l_2 $ are
    \[ \begin {array}{l|cc}
         \mbox {$ \downarrow $ evaluation at\dots} & l_1 & l_2 \\ \hline
         \mbox {$ x_1 $ (1 row)} & v^1 & 0 \\
         \mbox {points reached via $ E_1 $ from $ x_1 $ (2 rows each,
           except only 1 for $ x_2 $)}
           & 0 & 0 \\
         \mbox {points reached via $ E_2 $ from $ x_1 $ (2 rows each,
           except only 1 for $ x_2 $)}
           & v & v
       \end {array} \]
  We see that after subtracting the $ l_2 $-column from the $ l_1 $-column we
  again get one column with only one non-zero entry $ v^1 $. So for the
  determinant we get $ v^1 = (C_1 \cdot L_1)_{x_1} $ as a factor, dropping the
  corresponding row and column (which simply means forgetting the point $ x_1 $
  as in definition \ref {def-forget} \ref {def-forget-b}). Essentially the same
  argument holds if one of the adjacent edges --- say $ E_2 $ -- is unbounded:
  in this case there is only an $ l_1 $-column which has zeroes everywhere
  except in the one $ x_1 $-row where the entry is $ v^1 $.

  The same is of course true for $ x_2 $ and leads to a factor of $ (C_1 \cdot
  L_2)_{x_2} $.

  Next we consider again the contracted bounded edge $E$ at which we split the
  curve $C$ into the two parts $ C_1 $ and $ C_2 $. Choose one of its
  boundary points as root vertex $V$, say the one on the $ C_1 $ side. Denote
  the adjacent edges and their directions as in the following picture:

  \begin {center} \input {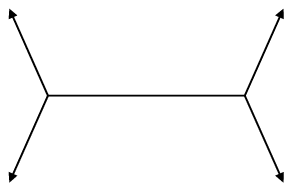} \end {center}

  If we set $ l_i = l(E_i) $ the matrix of $ \pi $ (of size $ 2n-4 $) reads
    \[ \begin {array}{c|l|ccc@{\;\,}c@{\;\,}c@{\;\,}cc}
         & & & \mbox {lengths in $ C_1 $} & & & & & 
           \mbox {lengths in $ C_2 $} \\
         & & \mbox {root} &
           \mbox {($ 2n_1-3 $ cols)} & l_1 & l_2 & l_3 & l_4 &
           \mbox {($ 2n_2+1 $ cols)} \\ \hline \hline
       (2n_1 \; \hphantom {+4} &
         \mbox {pts behind $ E_1 $} &
         I_2 & * & v & 0 & 0 & 0 & 0 \\
       \mbox {rows}) & \mbox {pts behind $ E_2 $} &
         I_2 & * & 0 & -v & 0 & 0 & 0 \\ \hline
       (2n_2+4 &
         \mbox {pts behind $ E_3 $} &
         I_2 & 0 & 0 & 0 & w & 0 & * \\
       \mbox {rows}) & \mbox {pts behind $ E_4 $} &
         I_2 & 0 & 0 & 0 & 0 & -w & *
       \end {array} \]
  where $ n_1 $ and $ n_2 $ are as in the proof of lemma \ref {set-reducible},
  $ I_2 $ is the $ 2 \times 2 $ unit matrix, and $*$ denotes arbitrary entries.
  Now add $v$ times the root columns to the $ l_2 $-column, subtract the $ l_1
  $-column from the $ l_2 $-column and the $ l_4 $-column from the $ l_3
  $-column to obtain the following matrix with the same determinant:
    \[ \begin {array}{c|l|ccc|c@{\;\,}c@{\;\,}cc}
         & & & \mbox {lengths in $ C_1 $} & & & & & 
           \mbox {lengths in $ C_2 $} \\
         & & \mbox {root} &
           \mbox {($ 2n_1-3 $ cols)} & l_1 & l_2 & l_3 & l_4 &
           \mbox {($ 2n_2+1 $ cols)} \\ \hline \hline
       (2n_1 \; \hphantom {+4} &
         \mbox {pts behind $ E_1 $} &
         I_2 & * & v & 0 & 0 & 0 & 0 \\
       \mbox {rows}) & \mbox {pts behind $ E_2 $} &
         I_2 & * & 0 & 0 & 0 & 0 & 0 \\ \hline
       (2n_2+4 &
         \mbox {pts behind $ E_3 $} &
         I_2 & 0 & 0 & v & w & 0 & * \\
       \mbox {rows}) & \mbox {pts behind $ E_4 $} &
         I_2 & 0 & 0 & v & w & -w & *
       \end {array} \]
  Note that this matrix has a block form with a zero block at the top right.
  Denote the top left block (of size $ 2n_1 $) by $ A_1 $ and the bottom
  right (of size $ 2n_2+4 $) by $ A_2 $, so that the multiplicity that we are
  looking for is $ | \det A_1 \cdot \det A_2 | $.

  The matrix $ A_1 $ is precisely the matrix for the evaluation map of $ C_1 $
  if we forget the marked point corresponding to $E$ and choose the other end
  point of $ E_2 $ as the root vertex. Hence $ |\det A_1| = \mult_{\ev}(C_1) $.
  In the same way the matrix for the evaluation map of $ C_2 $, if we again
  forget the marked point corresponding to $E$ and now choose the other end
  point of $ E_3 $ as the root vertex, is the matrix $ A_2' $ obtained from $
  A_2 $ by replacing $v$ and $w$ in the first two columns by the first and
  second unit vector, respectively. But $ A_2 $ is simply obtained from $
  A_2' $ by right multiplication with the matrix
    \[ \left( \begin {array}{ccc}
         v & w & 0 \\
         0 & 0 & I_{2n_2+2}
       \end {array} \right) \]
  which has determinant $ \det (v,w) $. So we conclude that
    \[ |\det A_2| = |\det (v,w)| \cdot |\det A_2'|
         = (C_1 \cdot C_2)_{P=Q} \cdot \mult_{\ev}(C_2). \]
  Collecting these results we now obtain the formula stated in the proposition.
\end {proof}

Of course there are completely analogous statements to lemma \ref
{set-reducible}, remark \ref {converse-reducible}, and proposition \ref
{comp-mult} if the $ \calM_4 $-coordinate of the curves in question is of type
(B) instead of type (A) (see example \ref {ex-m4}). Note however that there are
no curves of type \ref {set-reducible-a} in lemma \ref {set-reducible} in this
case since $ x_1 $ and $ x_3 $ would have to map to $ L_1 \cap p_3 $, which is
empty.

We can now collect our results to obtain the final theorem. The idea of this
final step is actually the same as in the case of complex curves.

\begin {theorem}[Kontsevich's formula] \label {kf}
  The numbers $ N_d $ of example \ref {ex-eval} and remark \ref {rem-eval}
  satisfy the recursion formula
    \[ N_d = \sum_{\substack {d_1+d_2=d \\ d_1,d_2>0}} \left(
         d_1^2 d_2^2 \binom {3d-4}{3d_1-2} - d_1^3 d_2 \binom {3d-4}{3d_1-1}
       \right) N_{d_1} N_{d_2} \]
  for $ d>1 $.
\end {theorem}

\begin {proof}
  We compute the degree of the map $ \pi $ of definition \ref {def-pi} at two
  different points. First consider a point $ \calP=(a,b,p_3,\dots,p_n,z) \in
  \RR^{2n-2} \times \calM_4 $ in $ \pi $-general position with $ \calM_4
  $-coordinate $z$ of type (A) (see example \ref {ex-m4}) with a very large
  length. We have to count the points in $ \pi^{-1}(\calP) $ with their
  respective $ \pi $-multiplicity. Starting with the curves of type \ref
  {set-reducible-a} in lemma \ref {set-reducible} we see by proposition \ref
  {comp-mult} that they simply count curves of degree $d$ through $ 3d-1 $
  points with their ordinary ($ \ev $-)multiplicity, so this simply gives us a
  contribution of $ N_d $. For the curves of type \ref {set-reducible-b} remark
  \ref {converse-reducible} tells us that we can as well count tuples $
  (C_1,C_2,x_1,\dots,x_n,P,Q) $, where
  \begin {enumerate}
  \item \label {kf-a}
    $ C_1 $ and $ C_2 $ are tropical curves of degrees $ d_1 $ and $ d_2 $
    with $ d_1+d_2=d $;
  \item \label {kf-b}
    $ x_1,x_2 $ are marked points on $ C_1 $ that map to $ L_1 $ and $ L_2
    $, respectively;
  \item \label {kf-c}
    $ x_3,x_4 $ are marked points on $ C_2 $ that map to $ p_3 $ and $ p_4
    $, respectively;
  \item \label {kf-d}
    $ x_5,\dots,x_n $ are marked points that map to $ p_5,\dots,p_n $ and
    of which exactly $ 3d_1-1 $ lie on $ C_1 $ and $ 3d_2-3 $ on $ C_2 $;
  \item \label {kf-e}
    $ P \in C_1 $ and $ Q \in C_2 $ are points with the same image in $
    \RR^2 $;
  \end {enumerate}
  where each such tuple has to be counted with the multiplicity computed in
  proposition \ref {comp-mult}.

  There are $ \binom {3d-4}{3d_1-1} $ choices to split up the points $
  x_5,\dots,x_n $ as in \ref {kf-d}. After fixing $ d_1 $ and $ d_2 $ we then
  have $ N_{d_1} \cdot N_{d_2} $ choices for $ C_1 $ and $ C_2 $ in \ref {kf-a}
  if we count each of them with their $ \ev $-multiplicity (which we have to do
  by proposition \ref {comp-mult}). By B\'ezout's theorem (see \cite {RST03}
  theorem 4.2) there are $ d_1 $ possibilities for $ x_1 $ in \ref {kf-b} ---
  namely the intersection points of $ C_1 $ with $ L_1 $ --- if we count each
  of them with its local intersection multiplicity $ (C_1 \cdot L_1)_{x_1} $ as
  required by proposition \ref {comp-mult}. In the same way there are again $
  d_1 $ choices for $ x_2 $ and $ d_1 \cdot d_2 $ choices for the glueing point
  $ P=Q $. (Note that we can apply B\'ezout's theorem without problems since we
  have seen in remark \ref {rem-reducible} that $ C_1 $ intersects $ L_1 $, $
  L_2 $, and $ C_2 $ in only finitely many points.)

  Altogether we see that the degree of $ \pi $ at $ \calP $ is
    \[ \deg_\pi (\calP) = N_d + \sum_{d_1+d_2=d}
         d_1^3 d_2 \binom {3d-4}{3d_1-1} N_{d_1} N_{d_2}. \]
  Repeating the same arguments for a point $ \calP' $ with $ \calM_4
  $-coordinate of type (B) as in example \ref {ex-m4} we get
    \[ \deg_\pi (\calP') = \sum_{d_1+d_2=d}
         d_1^2 d_2^2 \binom {3d-4}{3d_1-2} N_{d_1} N_{d_2}. \]
  Equating these two expressions by proposition \ref {np-indep} now gives the
  desired result.
\end {proof}




\providecommand{\bysame}{\leavevmode\hbox to3em{\hrulefill}\thinspace}
\providecommand{\MR}{\relax\ifhmode\unskip\space\fi MR }

\providecommand{\href}[2]{#2}

\end {document}